\documentclass[12pt]{amsart}
\usepackage[utf8]{inputenc}

\usepackage{subfiles}
\usepackage{hyperref}
\usepackage{cleveref}
\usepackage{amsmath}
\usepackage{multirow}
\usepackage{amssymb}
\usepackage{amsthm}
\usepackage{mathtools}
\usepackage{array}
\usepackage{thmtools}
\usepackage{thm-restate}
\usepackage[margin=1in]{geometry}
\usepackage{cleveref}
\usepackage{ytableau}
\usepackage{comment}
\usepackage{soul}
\usepackage{tikz}
\usetikzlibrary{patterns, matrix}
\tikzstyle{NE-lines}=[pattern=north east lines, pattern color=black!45]
\usetikzlibrary{positioning}
\usepackage{booktabs}

\usepackage{graphicx}
\usepackage{caption}
\usepackage{float}

\usepackage[backend=bibtex]{biblatex}
\renewbibmacro{in:}{}
\DeclareFieldFormat[article]{citetitle}{#1}
\DeclareFieldFormat[article]{title}{#1}
\DeclareFieldFormat{pages}{#1}
\DeclareUnicodeCharacter{2212}{-}

\numberwithin{equation}{section}

\DeclareMathOperator{\ord}{ord}

\theoremstyle{definition}
\newtheorem{theorem}{Theorem}[section]
\newtheorem*{theorem*}{Theorem}

\newtheorem*{example*}{Example}
\newtheorem{lemma}[theorem]{Lemma}
\newtheorem*{lemma*}{Lemma}
\newtheorem{corollary}[theorem]{Corollary}
\newtheorem*{corollary*}{Corollary}

\newtheorem*{definition*}{Definition}
\newtheorem{proposition}[theorem]{Proposition}
\newtheorem*{proposition*}{Proposition}

\newtheorem*{remark*}{Remark}
\newtheorem{conjecture}[theorem]{Conjecture}

\usepackage{ mathrsfs }

\title[The (123, 132)-Avoiding Stack Sort]{The Order of the (123, 132)-Avoiding Stack Sort}
\author{Owen Zhang}\address{\textsc{O. Zhang}, Yale University,
    New Haven, CT} \email{owen.zhang.oz37@yale.edu}

\begin{document}
\begin{abstract}
Let $s$ be West's deterministic stack-sorting map. A well-known result (West) is that any length $n$ permutation can be sorted with $n-1$ iterations of $s.$ In 2020, Defant introduced the notion of \textit{highly-sorted} permutations---permutations in $s^t(S_n)$ for $t \lessapprox n-1.$ In 2023, Choi and Choi extended this notion to generalized stack-sorting maps $s_{\sigma},$ where we relax the condition of becoming sorted to the analogous condition of becoming periodic with respect to $s_{\sigma}.$ In this work, we introduce the notion of \textit{minimally-sorted} permutations $\mathfrak{M}_n$ as an antithesis to Defant's highly-sorted permutations, and show that $\ord_{s_{123, 132}}(S_n) = 2 \lfloor \frac{n-1}{2} \rfloor,$ strengthening Berlow's 2021 classification of periodic points.
\end{abstract}

\maketitle

\section{Introduction} 
Knuth's \textit{Art of Computer Programming} \cite{KNUTH1968} first introduced the \textit{stack-sorting machine}, in which an input sequence is sorted with a single external stack structure. The elements of the sequence are passed left-to-right through the machine, with two possible operations at every state: \textit{push}, moving the next input element onto the stack, and \textit{pop}, removing the top element from the stack and appending it to the output. 

In 1990, West \cite{WEST1990} introduced a deterministic version of Knuth's stack-sorting machine as the \textit{stack-sorting map} $s$, insisting that the stack must always increase from top to bottom and employ a \textit{right-greedy} process: the \textit{push} operation is prioritized. Since then, various studies have been motivated by Knuth's original machine and West's deterministic $s$, including pop-stack-sorting \cite{ATKINSON1999, AVIS1981, DEFANT2022, PUDWELL2019, SMITH2014}, stack-sorting Coxeter groups \cite{DEFANT2022, DEFANT20222}, sigma-tau machines \cite{BAO2023, BARIL2021, BERLOW2021112571}, and stack-sorting of set-partitions \cite{DEFANTKRAVITZ2022, XIA2023}.

\begin{figure}[h]
\begin{center}
\begin{tikzpicture}[scale=0.6]
				\draw[thick] (0,0) -- (1,0) -- (1,-2) -- (2,-2) -- (2,0) -- (3,0);
				\node[fill = white, draw = white] at (2.5,.5) {2143};
				\node[fill = white, draw = white] at (1.5,-1.5) {};
				\node[fill = white, draw = white] at (1.5,-.8) {};
				\node[fill = white, draw = white] at (3.5,-1) {$\rightarrow$};
              \end{tikzpicture}
              \begin{tikzpicture}[scale=0.6]
				\draw[thick] (0,0) -- (1,0) -- (1,-2) -- (2,-2) -- (2,0) -- (3,0);
				\node[fill = white, draw = white] at (2.5,.5) {143};
				\node[fill = white, draw = white] at (1.5,-1.5) {2};
                \node[fill = white, draw = white] at (1.5,-.8) {};
				\node[fill = white, draw = white] at (.5,.5) {};
             \node[fill = white, draw = white] at (3.5,-1) {$\rightarrow$}; \end{tikzpicture}
         \begin{tikzpicture}[scale=0.6]
				\draw[thick] (0,0) -- (1,0) -- (1,-2) -- (2,-2) -- (2,0) -- (3,0);
				\node[fill = white, draw = white] at (2.5,.5) {43};
				\node[fill = white, draw = white] at (1.5,-1.5) {2};
				\node[fill = white, draw = white] at (1.5,-.8) {1};
				\node[fill = white, draw = white] at (3.5,-1) {$\rightarrow$};
              \end{tikzpicture} 
              \begin{tikzpicture}[scale=0.6]
				\draw[thick] (0,0) -- (1,0) -- (1,-2) -- (2,-2) -- (2,0) -- (3,0);
				\node[fill = white, draw = white] at (2.5,.5) {43};
				\node[fill = white, draw = white] at (1.5,-1.5) {2};
				\node[fill = white, draw = white] at (1.5,-.8) {};
                \node[fill = white, draw = white] at (1.5,0) {};
                \node[fill = white, draw = white] at (.5,.5) {1};
             	\node[fill = white, draw = white] at (3.5,-1) {$\rightarrow$}; \end{tikzpicture}
              \begin{tikzpicture}[scale=0.6]
				\draw[thick] (0,0) -- (1,0) -- (1,-2) -- (2,-2) -- (2,0) -- (3,0);
				\node[fill = white, draw = white] at (2.5,.5) {43};
				\node[fill = white, draw = white] at (1.5,-1.5) {};
                 \node[fill = white, draw = white] at (1.5,-.8) {};
				\node[fill = white, draw = white] at (.5,.5) {12};
             	\node[fill = white, draw = white] at (3.5,-1) {$\rightarrow$}; \end{tikzpicture}\begin{tikzpicture}[scale=0.6]
				\draw[thick] (0,0) -- (1,0) -- (1,-2) -- (2,-2) -- (2,0) -- (3,0);
				\node[fill = white, draw = white] at (2.5,.5) {3};
				\node[fill = white, draw = white] at (1.5,-1.5) {4};
                 \node[fill = white, draw = white] at (1.5,-.8) {};
                 \node[fill = white, draw = white] at (1.5,-.1) {};
				\node[fill = white, draw = white] at (.5,.5) {12};
             	\node[fill = white, draw = white] at (3.5,-1) {$\rightarrow$}; \end{tikzpicture}\begin{tikzpicture}[scale=0.6]
				\draw[thick] (0,0) -- (1,0) -- (1,-2) -- (2,-2) -- (2,0) -- (3,0);
				\node[fill = white, draw = white] at (2.5,.5) {};
				\node[fill = white, draw = white] at (1.5,-1.5) {4};
                 \node[fill = white, draw = white] at (1.5,-.8) {3};
				\node[fill = white, draw = white] at (.5,.5) {12};
             	\node[fill = white, draw = white] at (3.5,-1) {$\rightarrow$}; \end{tikzpicture}\begin{tikzpicture}[scale=0.6]
				\draw[thick] (0,0) -- (1,0) -- (1,-2) -- (2,-2) -- (2,0) -- (3,0);
				\node[fill = white, draw = white] at (2.5,.5) {};
				\node[fill = white, draw = white] at (1.5,-1.5) {4};
                 \node[fill = white, draw = white] at (1.5,-.8) {};
				\node[fill = white, draw = white] at (.5,.5) {123};
             	\node[fill = white, draw = white] at (3.5,-1) {$\rightarrow$}; \end{tikzpicture}
              \begin{tikzpicture}[scale=0.6]
				\draw[thick] (0,0) -- (1,0) -- (1,-2) -- (2,-2) -- (2,0) -- (3,0);
				\node[fill = white, draw = white] at (.5,.5) {1234};
\end{tikzpicture}
\end{center}
\caption{West's deterministic stack-sorting map $s$ on $\pi = 2143.$}
\label{Westmap}
\end{figure}
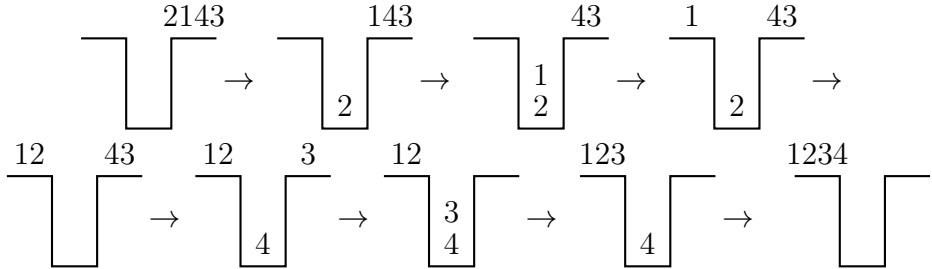

In his dissertation, West \cite{WEST1990} proved that $s^{n-1}(S_n)$ contains only the identity permutation, justifying repeated applications of $s$ as a correct and terminating sorting algorithm. A natural direction of study, then, is the characterization of \textit{$t$-stack-sortable} permutations---permutations $\pi$ such that $s^{t}(\pi)$ is sorted---for general $t \leq n-1$. Knuth \cite{KNUTH1968} answered the question for $t = 1,$ showing that $\pi$ is $1$-stack-sortable if and only if $\pi$ avoids subsequences of the pattern $231,$ enumerating the number of such permutations of length $n$ to be $\frac{1}{n+1} \binom{2n}{n},$ the $n^{\text{th}}$ Catalan number. In 1990, West \cite{WEST1990} characterized the $2$-stack-sortable permutations, proving that $\pi$ is $2$-stack-sortable if and only if $\pi$ avoids subsequences of the pattern $2341$ and the barred pattern $3\overline{5}241.$ He also conjectured that the number of such permutations of length $n$ is $\frac{2}{(n+1)(2n+1)} \binom{3n}{n},$ which was proven by Zeilberger \cite{ZEILBERGER19928592} two years later. West \cite{WEST1990,WEST1993} then searched for a polynomial $P(n)$ such that $3$-stack-sortable permutations could be enumerated by $\frac{1}{P(n)} \binom{4n}{n},$ but was unsuccessful for $\mathrm{deg}(P(n)) < 7.$ In 2012, Úlfarsson \cite{ULFARSSON2012} characterized $3$-stack-sortable permutations with ``decorated patterns," but only in 2021, did Defant \cite{DEFANT2021} discover a polynomial-time algorithm to enumerate $3$-stack-sortable permutations. 

In 2020, Defant \cite{DEFANT2020} first considered $t$-stack-sortable permutations to be duals of the $t$\textit{-sorted permutations} \cite{DEFANT2020511526}---permutations in the image of $s^t(S_n),$ a generalization of Bousquet-Mélou's definition \cite{BOUSQUET2000} of \textit{sorted.} Defant then defined a permutation $\pi \in S_n$ to be \textit{highly-sorted} if $\pi$ is $t$-sorted for some $t$ close to $n,$ proving that a $t$-sorted permutation can contain at most $\lfloor \frac{n-t}{2} \rfloor$ descents \cite{DEFANT2020}. 



The classical stack-sorting map $s$ has since been generalized to $s_{\sigma}$ \cite{CERBAI2020105230} for permutations $\sigma,$ where instead of insisting that the stack increases, we insist that the stack avoids top-to-bottom subsequences of the pattern $\sigma.$ In 2021, Berlow \cite{BERLOW2021112571} introduced the family of maps $s_T,$ where the stack must avoid top-to-bottom subsequences of every pattern in set $T$ (see \Cref{Berlowmap}). In 2023, Choi and Choi \cite{CHOI2023} generalized Defant's notion of highly-sorted permutations, defining $\pi$ to be \textit{highly-sorted} with respect to $s_{\sigma}$ if $\pi$ is in the image of $s^{t}_{\sigma}$ for some $t$ close to $\ord_{s_{\sigma}}(S_n),$ where $\ord_{s_{\sigma}}(P)$ is the smallest integer $k$ such that every element in $s_{\sigma}^{k}(P)$ is periodic under $s_{\sigma}$. We straightforwardly extend this definition to generalized maps $s_T.$

\begin{figure}[h]
\begin{center}
\footnotesize
\begin{tikzpicture}[scale=0.55]
				\draw[thick] (0,0) -- (1,0) -- (1,-2) -- (2,-2) -- (2,0) -- (3,0);
				\node[fill = white, draw = white] at (2.5,.5) {52431};
				\node[fill = white, draw = white] at (1.5,-1.5) {};
				\node[fill = white, draw = white] at (1.5,-.8) {};
				\node[fill = white, draw = white] at (3.5,-1) {$\rightarrow$};
              \end{tikzpicture}
              \begin{tikzpicture}[scale=0.55]
				\draw[thick] (0,0) -- (1,0) -- (1,-2) -- (2,-2) -- (2,0) -- (3,0);
				\node[fill = white, draw = white] at (2.5,.5) {2431};
				\node[fill = white, draw = white] at (1.5,-1.5) {5};
                \node[fill = white, draw = white] at (1.5,-.8) {};
				\node[fill = white, draw = white] at (.5,.5) {};
             \node[fill = white, draw = white] at (3.5,-1) {$\rightarrow$}; \end{tikzpicture}
         \begin{tikzpicture}[scale=0.55]
				\draw[thick] (0,0) -- (1,0) -- (1,-2) -- (2,-2) -- (2,0) -- (3,0);
				\node[fill = white, draw = white] at (2.5,.5) {431};
				\node[fill = white, draw = white] at (1.5,-1.5) {5};
				\node[fill = white, draw = white] at (1.5,-.8) {2};
				\node[fill = white, draw = white] at (3.5,-1) {$\rightarrow$};
              \end{tikzpicture} 
              \begin{tikzpicture}[scale=0.55]
				\draw[thick] (0,0) -- (1,0) -- (1,-2) -- (2,-2) -- (2,0) -- (3,0);
				\node[fill = white, draw = white] at (2.5,.5) {31};
				\node[fill = white, draw = white] at (1.5,-1.5) {5};
				\node[fill = white, draw = white] at (1.5,-.8) {2};
                \node[fill = white, draw = white] at (1.5, -.1) {4};
                \node[fill = white, draw = white] at (.5,.5) {};
             	\node[fill = white, draw = white] at (3.5,-1) {$\rightarrow$}; \end{tikzpicture}
              \begin{tikzpicture}[scale=0.55]
				\draw[thick] (0,0) -- (1,0) -- (1,-2) -- (2,-2) -- (2,0) -- (3,0);
				\node[fill = white, draw = white] at (2.5,.5) {31};
				\node[fill = white, draw = white] at (1.5,-1.5) {5};
                 \node[fill = white, draw = white] at (1.5,-.8) {2};
				\node[fill = white, draw = white] at (.5,.5) {4};
             	\node[fill = white, draw = white] at (3.5,-1) {$\rightarrow$}; \end{tikzpicture}\begin{tikzpicture}[scale=0.55]
				\draw[thick] (0,0) -- (1,0) -- (1,-2) -- (2,-2) -- (2,0) -- (3,0);
				\node[fill = white, draw = white] at (2.5,.5) {1};
				\node[fill = white, draw = white] at (1.5,-1.5) {5};
                 \node[fill = white, draw = white] at (1.5,-.8) {2};
                 \node[fill = white, draw = white] at (1.5,-.1) {3};
				\node[fill = white, draw = white] at (.5,.5) {4};
             	\node[fill = white, draw = white] at (3.5,-1) {$\rightarrow$}; \end{tikzpicture}\begin{tikzpicture}[scale=0.55]
				\draw[thick] (0,0) -- (1,0) -- (1,-2) -- (2,-2) -- (2,0) -- (3,0);
				\node[fill = white, draw = white] at (2.5,.5) {1};
				\node[fill = white, draw = white] at (1.5,-1.5) {5};
                 \node[fill = white, draw = white] at (1.5,-.8) {2};
				\node[fill = white, draw = white] at (.5,.5) {43};
             	\node[fill = white, draw = white] at (3.5,-1) {$\rightarrow$}; \end{tikzpicture}\begin{tikzpicture}[scale=0.55]
				\draw[thick] (0,0) -- (1,0) -- (1,-2) -- (2,-2) -- (2,0) -- (3,0);
				\node[fill = white, draw = white] at (2.5,.5) {1};
				\node[fill = white, draw = white] at (1.5,-1.5) {5};
                 \node[fill = white, draw = white] at (1.5,-.8) {};
				\node[fill = white, draw = white] at (.5,.5) {432};
             	\node[fill = white, draw = white] at (3.5,-1) {$\rightarrow$}; 
              \end{tikzpicture}
              \begin{tikzpicture}[scale=0.55]
				\draw[thick] (0,0) -- (1,0) -- (1,-2) -- (2,-2) -- (2,0) -- (3,0);
				\node[fill = white, draw = white] at (2.5,.5) {};
				\node[fill = white, draw = white] at (1.5,-1.5) {5};
                 \node[fill = white, draw = white] at (1.5,-.8) {1};
				\node[fill = white, draw = white] at (.5,.5) {432};
                \node[fill = white, draw = white] at (3.5,-1) {$\rightarrow$}; 
            \end{tikzpicture}
            \begin{tikzpicture}[scale=0.55]
				\draw[thick] (0,0) -- (1,0) -- (1,-2) -- (2,-2) -- (2,0) -- (3,0);
				\node[fill = white, draw = white] at (2.5,.5) {};
				\node[fill = white, draw = white] at (1.5,-1.5) {5};
                 \node[fill = white, draw = white] at (1.5,-.8) {};
				\node[fill = white, draw = white] at (.5,.5) {4321};
                \node[fill = white, draw = white] at (3.5,-1) {$\rightarrow$}; 
            \end{tikzpicture}
            \begin{tikzpicture}[scale=0.55]
				\draw[thick] (0,0) -- (1,0) -- (1,-2) -- (2,-2) -- (2,0) -- (3,0);
				\node[fill = white, draw = white] at (2.5,.5) {};
				\node[fill = white, draw = white] at (1.5,-1.5) {};
                 \node[fill = white, draw = white] at (1.5,-.8) {};
				\node[fill = white, draw = white] at (.5,.5) {43215};
            \end{tikzpicture}
\end{center}
\caption{The generalized stack-sorting map $s_{123, 132}$ on $\pi = 52431.$}
\label{Berlowmap}
\end{figure}

Recently, Choi, Gan, Li, and Zhu \cite{CHOI2024} studied set partitions that require the maximum number of sorts through an $aba$-avoiding stack. Similarly, we define a permutation $\pi$ to be \textit{minimally-sorted} with respect to $s_T$ if $\ord_{s_T}(S_n) = \ord_{s_T}(\{\pi\}),$ antithetical to Defant's notion of highly-sorted permutations. At the end of this work, we present two conjectures on $\mathfrak{M}_n,$ the minimally-sorted permutations with respect to $s_{123, 132}.$

In 2021, Berlow \cite{BERLOW2021112571} studied the periodic points of $s_{123, 132}.$ She defined a permutation $\pi$ of length $n$ to be \textit{half-decreasing} if the subsequence $\pi_{n-1} \pi_{n-3} \cdots \pi_{(3-(n \text{ mod } 2))}$ is the identity of length $\lfloor \frac{n-1}{2} \rfloor.$ In particular, being order-isomorphic to the identity is not sufficient. 

\begin{theorem}[Berlow \cite{BERLOW2021112571}] \label{berlowthm}
A permutation $\pi$ is periodic under $s_{123, 132}$ if and only if $\pi$ is half-decreasing.
\end{theorem}

Our main result is that we find the exact value of $\ord_{s_{123, 132}}(S_n),$ extending Berlow's work on periodic permutations. An analogous result for $s_{321, 312}$ follows directly from \Cref{main1}.

\begin{theorem} \label{main1}
For all positive integers $n,$ we have $\ord_{s_{123, 132}}(S_n) = 2 \lfloor \frac{n-1}{2} \rfloor.$
\end{theorem}

\section{Preliminaries}
We say that $a \in A$ is \textit{periodic} under $f: A \to B$ if there exists a positive integer $k$ such that $f^k(a) = a.$ For some ordered set $S,$ we use $S_i$ to denote the $i$th element of $S.$ 

Let $[n]$ denote $\{1, 2, \cdots, n\}$ for positive integers $n.$ A \textit{permutation}, written $\pi = \pi_1 \pi_2 \cdots \pi_n,$ is an ordering of distinct positive integers with \textit{length} $\mathrm{len}(\pi) = n.$ We say that $\pi_1, \pi_2, \cdots, \pi_n$ are the \textit{elements} of $\pi,$ and use $\pi_{[i:j]}$ to denote the subpermutation $\pi_i, \pi_{i+1}, \cdots, \pi_{j}.$ We define $\mathrm{ind}_{\pi}(i),$ the \textit{index of $i$} in $\pi$, to be $j,$ where $\pi_j = i.$ Let $S_n$ be the set of permutations with elements $[n].$ The \textit{reduction} of a permutation $\pi$ (equivalently, the \textit{standardization} \cite{DEFANT2020}), is the unique permutation $\mathrm{red}(\pi) \in S_n$ such that $\mathrm{red}(\pi)_i= j$ for $1 \leq i \leq n,$ where $\pi_i$ is the $j$th smallest number in $\{\pi_1, \pi_2, \cdots, \pi_n\}.$ Two permutations $\pi$ and $\sigma$ are \textit{order-isomorphic} if $\mathrm{red}(\pi) = \mathrm{red}(\sigma),$ and we write $\pi \cong \sigma.$ For instance, $\pi = 57816$ and $\sigma = 48917$ are order-isomorphic, since $\mathrm{red}(\pi) = \mathrm{red}(\sigma) = 24513.$ Given permutations $\pi$ and $\sigma,$ we say that $\pi$ \textit{contains} the pattern $\sigma$ if there exists a sequence of positive integers $a_1 < a_2 < \cdots < a_k$ such that $\pi' = \pi_{a_1} \pi_{a_2} \cdots \pi_{a_k} \cong \sigma.$ Otherwise, we say that $\pi$ \textit{avoids} $\sigma$ (equivalently, is $\sigma$\textit{-avoiding}). For instance, $\pi = 24513$ contains $\sigma = 132$ since $\pi_1 \pi_3 \pi_5 = 253 \cong \sigma,$ but avoids $\tau = 321.$ We use $\pi \cdot \tau$ to denote the \textit{concatenation} of $\pi$ and $\tau,$ and let $\mathrm{rev}(\pi)$ denote the \textit{reverse} of $\pi,$ namely $\pi_n \pi_{n-1} \cdots \pi_1.$

Next, an element $\pi_i$ of $\pi \in S_n$ is \textit{small} if $\pi_j \leq \lfloor \frac{n-1}{2} \rfloor.$ An element $\pi_i$ is a \textit{left-to-right minimum} (equivalently, \textit{ltr-min}) of $\pi$ if $\pi_i = \min(\pi_{[1:i]}).$ Additionally, we say that $\pi_i$ is a \textit{valley} if $\pi_i$ is a ltr-min, $\pi_{i+1}$ (if $i+1 \leq n$) is not a ltr-min, and $\pi_{i+2}$ (if $i+2 \leq n$) is a ltr-min. A consecutive subsequence of elements $\pi_{[i:i+j]}$ is a \textit{valley-block} $\overline{v}$ if $\pi_{i+j}$ is a valley and $\mathrm{red}(\pi_{[1:i+j]})_{[i:i+j]} = j+1, j, \cdots, 1.$ We say that the \textit{valley-boundary} of $\pi \in S_n,$ denoted $\mathfrak{B}(\pi),$ is the smallest index $i$ such that $\pi_{[i:n]} = \overline{v_1} \pi_{a_1} \overline{v_2} \pi_{a_2} \cdots \overline{v_j} \pi_{a_j}$ for valleys $\overline{v_1}, \cdots, \overline{v_j}$ and elements $\pi_{a_1}, \cdots, \pi_{a_j},$ and set $\mathfrak{B}(\pi) = n$ if no such index exists. The \textit{valley-region} of $\pi$ is $\pi_{[\mathfrak{B}(\pi):n]}.$ For instance, given $\pi = (11, 12, 7, 5, 8, 4, 3, 6, 2, 9, 1, 10),$ the elements $1, 2, 3,$ and $5$ are valleys and the sets $(7, 5), (4, 3), (2), (1)$ form $4$ valley-blocks in $\pi.$ Finally, $\mathfrak{B}(\pi) = 3,$ since $\pi_{[3:n]} = \overline{7,5}, 8, \overline{4,3}, 6, \overline{2}, 9, \overline{1}, 10.$

We conclude by noting that permutation indices will be considered modulo $n$ for the duration of this paper. In particular, let $\pi_i := \pi_j,$ where $j$ is the unique element of $[n]$ such that $i \equiv j \pmod{n}.$

\section{Proof of the Main Result}
We preface this section with two propositions, immediate from the preliminaries.

\begin{proposition} \label{prop}
Given $\sigma, \tau \in S_3,$ it holds that $(s_{\sigma, \tau}(\pi))_n = \pi_1$ for all $\pi \in S_n$ and $n \geq 1.$
\end{proposition}

\begin{proposition} \label{valleys}
Let $\overline{v_1}, \cdots, \overline{v_i}$ be the valley-blocks of $\pi$ from left to right, and let $\mathrm{len}(v_j) = l_j$ for all $j.$ Then, the permutation $\overline{v_1} \cdot \overline{v_2} \cdot \cdots \overline{v_i}$ is the reverse of the identity of length $\sum l_i.$
\end{proposition}

We now begin the proof of \Cref{main1} with several auxiliary lemmas that demonstrate the monovariant movement of valley-blocks under $s_{123, 132}.$ 

\begin{lemma} \label{fixedlocal} 
For any $\pi \in S_n$ and ltr-min $\pi_i$ with $i > 1,$ let $j \leq n$ be the largest index such that $\pi_i = \min(\pi_{[1:j]}).$ It holds that $s_{123, 132}(\pi)_{j-1} = \pi_i.$
\end{lemma}

\begin{proof}
Since $\pi_i$ is a ltr-min, just before $\pi_i$ enters the stack, $\pi_1$ must be the only element in the stack. After the elements $\pi_{[i+1:j]}$ have all entered the stack, $\pi_i$ and $\pi_1$ necessarily remain in the stack since $\pi_{i+1}, \cdots, \pi_j > \pi_i.$ Additionally, since $\pi_{j+1} < \pi_i,$ just before $\pi_{j+1}$ enters the stack, $\pi_j$ must exit the stack. At this moment, the $j-1$ elements $\pi_2, \pi_3, \cdots, \pi_{j}$ have been the only elements to exit the stack, with $\pi_i$ being the last, so $s_{123, 132}(\pi)_{j-1} = \pi_i.$
\end{proof}

\begin{lemma} \label{leftshift} 
Given a valley-block $\bar{v} = \pi_{[i:i+j]}$ of $\pi,$ we have $s_{123, 132}(\pi)_{i+j} = \pi_{i+j}$ and $s_{123, 132}(\pi)_{k-1} = \pi_k$ for $i \leq k < i+j.$
\end{lemma}

\begin{proof}
Just before $\pi_i$ enters the stack, $\pi_1$ must be the only element in the stack. Since $\overline{v}$ consists of the $j+1$ smallest elements of $\pi_{[1:i+j]}$ in descending-order, just before any element of $\overline{v}$ enters the stack, the previous element must exit. Hence, $k-2$ elements exit before $\pi_{k}$ for $i \leq k < i+j,$ and thus $s_{123, 132}(\pi)_{k-1} = \pi_k.$ Finally, by \Cref{fixedlocal}, $\pi_{i+j}$ is a fixed point.
\end{proof}

Next, we show that $s_{123, 132}$ preserves the elements in the valley-region of $\pi.$

\begin{lemma} \label{invariantexpansion2}
Suppose $\pi_{[i:j]}$ and $\pi_{[j+2:k]}$ are two valley-blocks of $\pi.$ Then, $s_{123, 132}(\pi)_{j-1} = \pi_{j+1}.$
\end{lemma}

\begin{proof}
Right before $\pi_{j}$ enters the stack, the only element remaining must be $\pi_1.$ Now, since $\pi_{j+1} > \pi_{j},$ the stack will read $\pi_{j+1} \pi_{j} \pi_1$ top to bottom just after $\pi_{j+1}$ enters. Finally, since $\pi_{j+2}$ is also a ltr-min, just before it enters, $\pi_{j+1}$ and $\pi_{j}$ must have left the stack. Hence, every element in $\pi_{[1:j]}$ exits the stack before $\pi_{j+1}$ except $\pi_1$ and $\pi_{j},$ yielding $s(\pi)_{j-1} = \pi_{j+1}.$
\end{proof}

\begin{lemma} \label{invariantexpansion1} 
If $\pi_i$ is in the valley-region of $\pi,$ then $\pi_i$ is also in the valley-region of $s_{123, 132}(\pi).$
\end{lemma} 
\begin{proof}
Let $\pi_{[\mathfrak{B}(\pi):n]} = \overline{v_1} \pi_{a_1} \overline{v_2} \cdots \overline{v_j} \pi_{a_j},$ the valley-region of $\pi,$ and let $\mathrm{len}(\overline{v_i}) = l_i$ for $1 \leq i \leq j.$ Then, by \Cref{fixedlocal} and \Cref{leftshift}, we have that $s_{123, 132}(\pi)$ ends with the suffix $({v_1}_{[1:l_1-1]}) \cdot \pi_{b_1} \cdot ({v_1}_{[l_1]} \cdot {v_2}_{[1:l_2-1]}) \cdot \pi_{b_2} \cdot ({v_2}_{[l_2]} \cdot {v_3}_{[1:l_3-1]}) \cdots ({v_{j-1}}_{[l_{j-1}]} \cdot {v_j}_{[1]}) \cdot \pi_{b_{j-1}} \cdot ({v_j}_{[l_j]}) \cdot \pi_{b_j}$ for some elements $\pi_{b_1}, \pi_{b_2}, \cdots, \pi_{b_j}.$ By \Cref{valleys}, this suffix is of the form $\overline{w_1} \pi_{b_{c_1}} \overline{w_2} \pi_{b_{c_2}} \cdots \overline{w_k} \pi_{b_{c_k}},$ where $\pi_{{b}_{c_1}}, \cdots, \pi_{b_{c_k}}$ are the elements of $\{\pi_{b_1}, \cdots, \pi_{b_j}\}$ that are not ltr-mins. Hence, this suffix is fully contained in the valley-region of $s_{123, 132}(\pi).$ However, it also contains all the elements in valley-blocks in $\pi_{[\mathfrak{B}(\pi):n]},$ and all the elements in between valley-blocks in $\pi_{[\mathfrak{B}(\pi):n]}$ by \Cref{invariantexpansion2}, which fully encompass all of elements in the valley-block, finishing the proof.
\end{proof}

\begin{lemma} \label{smallinvariance}
Let $\pi_i = \min(\pi_{[1:\mathfrak{B}(\pi)-1]}).$ If $\pi_i$ is small, then $\pi_i$ is in the valley-region of $s_{123, 132}(\pi).$
\end{lemma}

\begin{proof}
If $i = 1,$ then the claim follows from \Cref{prop}. Otherwise, just before $\pi_i$ enters the stack, $\pi_1$ must be the only element remaining in the stack, since $\pi_i$ is a ltr-minimum. Then, after $\pi_{i+1}, \cdots, \pi_{\mathfrak{B}(\pi)-1}$ have all entered the stack, $\pi_i$ will remain in the stack. However, when $\pi_{\mathfrak{B}(\pi)}$ enters the stack, $\pi_i$ will necessarily leave, since $\pi_{\mathfrak{B}}$ is part of a valley-block to the right of $\pi_i,$ so $\pi_{\mathfrak{B}(\pi)} < \pi_i.$ Thus, since every other element in $\pi_1, \cdots, \pi_{\mathfrak{B}(\pi)-1}$ was popped out before $\pi_i,$ except for $\pi_1,$ we have $s(\pi)_{\mathfrak{B}(\pi)-2} = \pi_i.$ However, since $\pi_i= \min(\pi_1, \cdots, \pi_{\mathfrak{B}(\pi)-1}),$ the proof of \Cref{invariantexpansion1} shows that $\pi_i$ is in the valley-region of $s_{123, 132}(\pi).$
\end{proof}

By \Cref{invariantexpansion1}, elements never leave the valley-region, and by \Cref{smallinvariance}, a small element is always added to the valley-region every iteration, implying the following result.

\begin{corollary} \label{msorts}
For any $\pi \in s_{123, 132}^{\lfloor \frac{n-1}{2} \rfloor}(S_n),$ it holds that $i \geq \mathfrak{B}(\pi)$ for all small elements $\pi_i.$
\end{corollary}

\Cref{msorts} gives a characterization of the $\lfloor \frac{n-1}{2} \rfloor$-sorted permutations under a $s_{123, 132}$ map. We continue by showing that these permutations become periodic with at most $\lfloor \frac{n-1}{2} \rfloor$ further passes. 

\begin{lemma} \label{barrier}
For $\pi \in S_n$ and small element $i,$ if $\pi_{n-2i+2} = i$ and $i$ is in the valley-region of $\pi,$ then $s_{123, 132}(\pi)_{n-2i+1} = i.$
\end{lemma}

\begin{proof}
Suppose for the sake of contradiction that $\pi_{n-2i+2}$ is directly in between two valley-blocks, so that $\pi_{[j:n-2i+1]}$ is a valley-block for some $j \leq n-2i.$ By definition, $\pi_{n-2i+1}$ is a valley, and by \Cref{fixedlocal}, $s^{k}(\pi)_{n-2i+1} = \pi_{n-2i+1}$ for all $k.$ But this contradicts \Cref{berlowthm}, since we have $s^{k}(\pi)_{n-2i+1} \neq \pi_{n-2i+2} = i.$ Now, suppose that $\pi_{n-2i+2}$ is itself a valley. This similarly contradicts \Cref{berlowthm}, since we have $s^{k}(\pi)_{n-2i+2} = \pi_{n-2i+2}$ for all $k$ by \Cref{fixedlocal}.

Since $\pi_{n-2i+2}$ is in the valley-region of $\pi,$ the only remaining possibility is that $\pi_{n-2i+2}$ is part of a valley-block but not a valley. Hence, by \Cref{invariantexpansion1}, we have $s_{123, 132}(\pi)_{n-2i+1} = i,$ as desired.
\end{proof}

\begin{lemma} \label{2msorts}
For all positive integers $i \leq \lfloor \frac{n-1}{2} \rfloor$ and $\pi \in S_n,$ the permutation \newline $\sigma_{n-1} \sigma_{n-3} \cdots \sigma_{n-2i+1}$ is the identity of length $i,$ where $\sigma = s^{i+\lfloor \frac{n-1}{2} \rfloor}_{123, 132}(\pi).$
\end{lemma}

\begin{proof}
We induct on $i.$ The base case $i = 1$ is immediate---in particular, $s^{\lfloor \frac{n-1}{2} \rfloor}_{123, 132}(\pi)_{n-1} = 1,$ which becomes a fixed element by \Cref{fixedlocal}, since otherwise $\mathfrak{B}(\pi) = n$ which contradicts \Cref{msorts}.

Now suppose that for some $1 < j \leq \lfloor \frac{n-1}{2} \rfloor,$ it holds that for all $\pi$ and $i < j,$ the permutation $\sigma_{n-1} \sigma_{n-3} \cdots \sigma_{n-2i+1}$ is the identity of length $i,$ where $\sigma = s_{123, 132}^{i+ \lfloor \frac{n-1}{2} \rfloor}(\pi).$ First, we note that by \Cref{leftshift} and \Cref{invariantexpansion2}, if an element $\pi_i$ is in the valley-region of $\pi,$ we have $s_{123, 132}(\pi)_{x} = \pi_i$ for some $x \in \{i-2, i-1, i\}.$ Next, consider some $\pi \in S_n,$ and let $\mathfrak{Z} = \{\mathrm{ind}_{s^{k}_{123, 132}(\pi)}(j) \mid \lfloor \frac{n-1}{2} \rfloor \leq k \leq \lfloor \frac{n-1}{2} \rfloor + j\}.$ By \Cref{fixedlocal}, if $\mathfrak{Z}_{l} - \mathfrak{Z}_{l+1} = 0$ for some $l \leq j,$ we must have $s^{\lfloor \frac{n-1}{2} \rfloor + l}_{123, 132}(\pi)_{n-2j+1} = j,$ or equivalently $\mathfrak{Z}_l \leq n-2j+1.$ Similarly, if $\mathfrak{Z}_{l} - \mathfrak{Z}_{l+1} = 1,$ we have by \Cref{leftshift} and \Cref{invariantexpansion2} that the element $j$ must be in a valley-block (but not a valley) of $s_{123, 132}^{\lfloor \frac{n-1}{2} \rfloor + l-1}(\pi),$ so by the inductive hypothesis, $\mathfrak{Z}_{l} \leq n-j-l+2.$ Otherwise, $\mathfrak{Z}_{l} - \mathfrak{Z}_{l+1} = 2,$ so we conclude recursively that $\mathfrak{Z}_{j+1} \leq n-2j+1.$ But combining \Cref{fixedlocal}, \Cref{barrier}, and the fact that $\mathfrak{Z}_{l} - \mathfrak{Z}_{l+1} \leq 2$ for all $l,$ we derive $\mathfrak{Z}_{j+1} = n-2j+1,$ or equivalently $s_{123, 132}^{\lfloor \frac{n-1}{2} \rfloor + j}(\pi)_{n-2j+1} = j.$ Hence, for all $\pi$ and $i < j+1,$ the permutation $\sigma_{n-1} \sigma_{n-3} \cdots \sigma_{n-2i+1}$ is the identity of length $i,$ where $\sigma = s_{123, 132}^{i + \lfloor \frac{n-1}{2} \rfloor}(\pi),$ completing the induction.

\end{proof}

In particular, any $\pi \in s_{123, 132}^{2 \lfloor \frac{n-1}{2} \rfloor}(S_n)$ is half-decreasing, which implies the following by \Cref{berlowthm}.

\begin{corollary} \label{upper}
For all positive integers $n,$ we have $\mathrm{ord}_{s_{123, 132}}(S_n) \leq 2 \lfloor \frac{n-1}{2} \rfloor.$
\end{corollary}

Finally, we present a family of minimally-sorted permutations to show that precisely $2 \lfloor \frac{n-1}{2} \rfloor$ iterations are required to sort all of $S_n.$ Define \[\gamma_n =  \left(\frac{n+1}{2}, 2, 3, \cdots, \frac{n-1}{2}, \frac{n+3}{2}, \cdots, n-2, 1, n-1, n\right)\] for odd $n \geq 5$ and $\gamma_n = \gamma_{n-1} \cdot n$ for even $n \geq 6.$ It is immediate that $\mathrm{ord}_{s_{123, 132}}([n]) = 2 \lfloor \frac{n-1}{2} \rfloor $ for $n \leq 4.$ Hence, we consider $n \geq 5.$ Let $\delta_n$ denote the permutation $\mathrm{rev}((\gamma_n)_{[2:n-3]})$ when $n$ is odd and $\mathrm{rev}((\gamma_n)_{[2:n-4]})$ when $n$ is even.

\begin{table}[h]
\centering
\begin{tabular}{ |c|c|} 
 \hline
 $n$ & $\gamma_n$ \\
 \hline
 $5$ & $(3,2,1,4,5)$ \\
 \hline
 $6$ & $(3,2,1,4,5,6)$ \\
 \hline
 $7$ & $(4,2,3,5,1,6,7)$ \\
 \hline
 $8$ & $(4,2,3,5,1,6,7,8)$ \\
 \hline
 $9$ & $(5,2,3,4,6,7,1,8,9)$ \\
 \hline
\end{tabular}
\caption{The first few $\gamma_n$ for $n \geq 5.$}
\label{gamma}
\end{table}

\begin{lemma} \label{equality}
For positive integers $n \geq 5$ and $k \leq \lfloor \frac{n-1}{2} \rfloor,$ we have $s^k_{123, 132}(\gamma_n)_{[1:n-2k-2]} = (\delta_n)_{[k:n-k-3]}$ for odd $n$ and $s^k_{123, 132}(\gamma_n)_{[1:n-2k-3]} = (\delta_n)_{[k:n-k-4]}$ for even $n.$ Furthermore, $\zeta_{n-1} {\zeta}_{n-3} \cdots {\zeta}_{n-2k+1}$ is the identity permutation of length $k,$ where $\zeta = s^k_{123, 132}(\gamma_n).$ 
\end{lemma}
\begin{proof}
We induct on $k.$ For brevity, we will prove the lemma for when $n$ is odd---the proof for even $n$ is directly analogous. For the base case $k = 1,$ we have $s_{123, 132}(\gamma_n)_{n} = (\gamma_n)_{1} = \frac{n+1}{2}$ by \Cref{prop}. Since $(\gamma_n)_{[2:n-3]}$ is strictly increasing, these elements are popped out in reverse order just before $1$ enters the stack. Hence, $s_{123, 132}(\gamma_n)_{[1:n-4]} = \delta_n = (\delta_n)_{[1:n-4]}.$ Finally, $s_{123, 132}(\gamma_n)_{n-1} = 1$ by \Cref{fixedlocal}, completing the base case.

Next, suppose $s^k_{123, 132}(\gamma_n)_{[1:n-2k-2]} = (\delta_n)_{[k:n-k-3]}$ for some $k$ and $\zeta_{n-1} \zeta_{n-3} \cdots \zeta_{n-2k+1}$ is the identity of length $k$ where $\zeta = s^{k}_{123, 132}(\gamma_n).$ By \Cref{prop}, we have $s^{k+1}_{123, 132}(\gamma_n)_n = s^k_{123, 132}(\gamma_n)_1,$ and since $s^k_{123, 132}(\gamma_n)_{[1:n-2k-2]}$ is strictly decreasing, it follows that these elements will exit the stack in the same order, giving $s^{k+1}_{123, 132}(\gamma_n)_{[1:n-2k-4]} = (\delta_n)_{[k+1:n-k-4]}$ by the inductive hypothesis. Finally, by \Cref{fixedlocal}, we have $s^{k+1}_{123, 132}(\gamma_n)_{n-2k-1} = k+1,$ completing the induction.
\end{proof}

\begin{lemma} \label{lower}
For all positive integers $n,$ we have $\mathrm{ord}_{s_{123, 132}}(S_n) \geq 2 \lfloor \frac{n-1}{2} \rfloor.$
\end{lemma}

\begin{proof}
It follows from \Cref{equality} that $s^{\lfloor \frac{n-1}{2} \rfloor-1}_{123, 132}(\gamma_n)_1 = \lfloor \frac{n-1}{2} \rfloor.$ By \Cref{prop} and \Cref{invariantexpansion2}, we have $\mathrm{ind}_{s^k_{123, 132}(\gamma_n)}(\lfloor \frac{n-1}{2} \rfloor) = n-2(k-\lfloor \frac{n-1}{2} \rfloor)$ for $k \geq \lfloor \frac{n-1}{2} \rfloor.$ Hence, $k = 2 \lfloor \frac{n-1}{2} \rfloor$ is the minimal $k$ such that $s^k_{123, 132}(\gamma_n)$ is half-decreasing, giving us the desired bound.
\end{proof}

Finally, we conclude that exactly $2 \lfloor \frac{n-1}{2} \rfloor$ iterations are required to sort $S_n.$

\begin{proof}[Proof of \Cref{main1}]
\Cref{upper} and \Cref{lower} directly imply $\mathrm{ord}_{s_{123, 132}}(S_n) = 2 \lfloor \frac{n-1}{2} \rfloor.$
\end{proof}

\section{Future Directions}
To study Defant's notion of highly-sorted permutations and our newly-introduced notion of minimally-sorted permutations, characterizing the periodic permutations under generalized stack-sorting maps is a prerequisite. We state a conjecture on the periodic points of other $s_{\sigma, \tau}$ stack-sorting maps for three pairs of $(\sigma, \tau)$, and restate a conjecture from Berlow.

\begin{conjecture}
For $(\sigma, \tau) = (123, 213), (132, 312), (231, 321),$ the map $s_{\sigma, \tau}$ is a bijection from $S_n$ to itself, and all permutations are periodic.
\end{conjecture}

\begin{conjecture} [Berlow \cite{BERLOW2021112571}]
For $(\sigma, \tau) = (213, 231), (132, 213), (231, 312),$ the only periodic points of $s_{\sigma, \tau}$ are the identity permutation and its inverse.
\end{conjecture}

Recall that $\mathfrak{M}_n$ is the set of minimally-sorted permutations under $s_{123, 132}.$ We conjecture several properties of elements in $\mathfrak{M}_n.$ However, these conditions are not sufficient for $n \geq 7.$

\begin{conjecture} \label{necessary}
For $\pi \in \mathfrak{M}_n,$ the following conditions hold true:
\begin{itemize}
    \item $\pi_1 \geq \lfloor \frac{n+1}{2} \rfloor.$ 
    \item \textbf{For odd $n$}: $\pi_{n-2} = 1$ and $\pi_{n-1}, \pi_n \geq \lfloor \frac{n+1}{2} \rfloor.$
    \item \textbf{For even $n$}: $\pi_{n-3} = 1$ and $\pi_{n-2}, \pi_{n-1}, \pi_n \geq \lfloor \frac{n+1}{2} \rfloor.$
\end{itemize}
\end{conjecture}

Next, an enumerative conjecture on $\mathfrak{M}_n,$ computationally verified for $n \leq 6.$

\begin{conjecture} \label{main2}
 For all positive integers $n,$ we have $|\mathfrak{M}_{2n}| = (n+1)|\mathfrak{M}_{2n-1}|.$
\end{conjecture}

Finally, we conclude with an enumerative conjecture on $\mathrm{Sort}_{t, n}(123, 132),$ the set of length $n$ permutations that are $t$-stack-sortable under $s_{123, 132}.$

\begin{conjecture} \label{enum}
For any positive integer $t$ and $n \geq 2t+1,$ we have:
\begin{itemize}
\item $|\mathrm{Sort}_{t, n}(123, 132)| = \frac{n+3}{2}|\mathrm{Sort}_{t, n-2}(123, 132)|$ if $n$ is odd.
\item $|\mathrm{Sort}_{t, n}(123, 132)| = \frac{n+4}{2}|\mathrm{Sort}_{t, n-2}(123, 132)|$ if $n$ is even.
\end{itemize}
\end{conjecture}

\section*{Acknowledgements}
The author thanks Yunseo Choi for suggesting the problem and giving helpful advice on the presentation of the paper.

\nocite{*}
\printbibliography
\end{document}